\documentclass[11pt]{article}
\usepackage[intlimits]{amsmath}
\usepackage{amsfonts, amsthm}
\usepackage{array}
\usepackage{amscd}
\usepackage{enumerate}
\usepackage{multicol}

\newtheorem{definition}{Definition}
\newtheorem{theorem}{Theorem}

\newcounter{ex}
\newenvironment{example}{\par\noindent\textbf{Example \stepcounter{ex}\arabic{ex}.}}{\bigskip}
\input{tcilatex}
\begin{document}

\title{On the existence of bounded solutions for nonlinear second order
neutral difference equations}
\date{}
\author{Marek Galewski\thanks{%
Lodz University of Technology, Poland, email: marek.galewski@p.lodz.pl} \and %
Robert Jankowski\thanks{%
University of Bialystok and Lodz University of Technology, Poland, email:
rjjankowski@math.uwb.edu.pl} \and Magdalena Nockowska--Rosiak\thanks{%
Lodz University of Technology, Poland, email: magdalena.nockowska@p.lodz.pl}
\and Ewa Schmeidel\thanks{%
University of Bialystok, Poland, email: eschmeidel@math.uwb.edu.pl} }
\maketitle

\begin{abstract}
Using the techniques connected with the measure of noncompactness we
\noindent investigate the neutral difference equation of the following form 
\begin{equation*}
\Delta \left( r_{n}\left( \Delta \left( x_{n}+p_{n}x_{n-k}\right) \right)
^{\gamma }\right) +q_{n}x_{n}^{\alpha }+a_{n}f(x_{n+1})=0,
\end{equation*}%
where $x:{\mathbb{N}}_{k}\rightarrow {\mathbb{R}}$, $a,p,q:{\mathbb{N}}%
_{0}\rightarrow {\mathbb{R}}$, $r:{\mathbb{N}}_{0}\rightarrow {\mathbb{R}}%
\setminus \{0\}$, $f\colon {\mathbb{R}}\rightarrow {\mathbb{R}}$ is
continuous and $k$ is a given positive integer, $\alpha \geq 1$ is a ratio
of positive integers with odd denominator, and $\gamma \leq 1$ is ratio of
odd positive integers; ${\mathbb{N}}_{k}:=\left\{ k,k+1,\dots \right\} $.
Sufficient conditions for the existence of a bounded solution are obtained.
Also a special type of stability and asymptotic stability are studied. Some
earlier results are generalized. We note that the solution which we obtain
does not directly correspond to a fixed point of a certain continuous
operator since it is partially iterated. The method which we develop allows
for considering through  techniques connected with the measure of
noncompactness also difference equations with memory. \newline
\end{abstract}

\section{Introduction}

As it is well known difference equations serve as mathematical models in
diverse areas, such as economy, biology, physics, mechanics, computer
science, finance, see for example \cite{bib1B}, \cite{bib2B}. One of such
models is the Emden--Fowler equation which originated in the gaseous
dynamics in astrophysics and further was used in the study of fluid
mechanics, relativistic mechanics, nuclear physics and in the study of
chemically reacting systems, see \cite{wong}. For the reader's convenience,
we note that the background for difference equations theory can be found in
numerous well--known monographs: Agarwal~\cite{bib1B}, Agarwal, Bohner,
Grace and O'Regan \cite{ABGO}, Agarwal and Wong \cite{AW} Elaydi~\cite{bib2B}%
, Kelley and Peterson~\cite{bib4B}, and Koci\'{c} and Ladas~\cite{bib3B}.

In the present paper we study using techniques connected with the measure of
noncompactness the existence of a bounded solution and some type of its
asymptotic behavior to a nonlinear second order difference equation, which
can be viewed as a generalization of a discrete Emeden--Fowler equation or
else it can be viewed as a second order difference equation with memory.
This makes a problem which we consider different from those already
investigated via techniques of measure of noncompactness, see for example 
\cite{SchmeildelZbaszyniakCAMW} since we do not expect a direct connection
between a fixed point of a suitably defined operator investigated on a
non-reflexive space $l^{\infty }$ a solution to the problem under
consideration. Indeed, this is the case. What we obtain is that starting
from some index which we define the solution is taken from the fixed point
while the previous terms have to be iterated. This also makes the definition
of the operator different from this which would be used had the problem been
without dependence on previous terms. It seems that the method which we
sketch here would prove applicable for several other problems. We also note
that due to the type of space which we use, namely $l^{\infty }$ we cannot
apply standard fixed point techniques such us Banach Theorem or Schauder
Theorem and related results. We expect that our method would apply for
systems of difference equations. However, what we cannot obtain here is the
asymptotic stability of the solution. This no surprise since solution which
we get is not is a fixed point of some mapping. We will use axiomatically
defined measures of noncompactness as presented in paper~\cite{bib1} by Bana%
\'{s} and Rzepka.

The problem we consider is as follows 
\begin{equation}
\Delta \left( r_{n}\left( \Delta \left( x_{n}+p_{n}x_{n-k}\right) \right)
^{\gamma }\right) +q_{n}x_{n}^{\alpha }+a_{n}f(x_{n+1})=0.  \label{e0}
\end{equation}%
where $\alpha \geq 1$ is a ratio of positive integers with odd denominator, $%
\gamma \leq 1$ is ratio of odd positive integers, $x:{\mathbb{N}}%
_{k}\rightarrow {\mathbb{R}}$ $a,p,q:{\mathbb{N}}_{0}\rightarrow {\mathbb{R}}
$, $r:{\mathbb{N}}_{0}\rightarrow {\mathbb{R}}\setminus \{0\}$, and $f:{%
\mathbb{R}}\rightarrow {\mathbb{R}}$ is a locally Lipschitz function with no
further growth assumptions. Here ${\mathbb{N}}_{0}:=\left\{ 0,1,2,\dots
\right\} $, ${\mathbb{N}}_{k}:=\left\{ k,k+1,\dots \right\} $ where $k$ is a
given positive integer, and ${\mathbb{R}}$ is a set of all real numbers. By
a solution of equation~\eqref{e0} we mean a sequence $x:{\mathbb{N}}%
_{k}\rightarrow {\mathbb{R}}$ which satisfies~\eqref{e0} for every $n\in {%
\mathbb{N}}_{k}$.

There has been an interest of many authors to study properties of solutions
of the second--order neutral difference equations attract attention; see the
papers \cite{DGJ2002}, \cite{bibGL}--\cite{JS2}, \cite{bibLG}--\cite{bibmm2}%
, \cite{S1}--\cite{S3}, \cite{TKP1}--\cite{TKP2} and the references therein.
The interesting oscillatory results for first order and even order neutral
difference equations can be found in \cite{LadasQian}, \cite{MM2004} and 
\cite{MM2009}.

\section{Preliminaries}

Let $(E,\left\Vert \cdot \right\Vert )$ be an infinite--dimensional Banach
space. If $X$ is a subset of $E$, then $\bar{X}$, $ConvX$ denote the closure
and the convex closure of $X$, respectively. Moreover, we denote by ${%
\mathcal{M}}_{E}$ the family of all nonempty and bounded subsets of $E$ and
by ${\mathcal{N}}_{E}$ the subfamily consisting of all relatively compact
sets.

\begin{definition}
A mapping $\mu\colon{\mathcal{M}}_E \to [0,\infty)$ is called a measure of
noncompactness in $E$ if it satisfies the following conditions:

\begin{description}
\item[$1^0$] $\ker \mu = \left\{X \in {\mathcal{M}}_E \colon \mu(X)=0
\right\}\neq \emptyset \mbox{  and  } \ker \mu \subset {\mathcal{N}}_E, $

\item[$2^0$] $X \subset Y \Rightarrow \mu(X) \leq \mu(Y), $

\item[$3^0$] $\mu(\bar{X})=\mu(X)=\mu(Conv \,\, X), $

\item[$4^0$] $\mu( c X +(1- c)Y ) \leq c \mu(X)+(1-c) \mu(Y) \mbox{  for  }
0 \leq c \leq 1, $

\item[$5^0$] If $X_n \in {\mathcal{M}}_E, \,\,\, X_{n+1}\subset X_n, \,\,\,
X_n=\bar{X_n} \mbox{  for  }n=1,2,3,\dots $\newline
and $\lim\limits_{n\to\infty} \mu(X_n)=0 \mbox{  then  } \bigcap
\limits_{n=1}^{\infty} X_n\neq \emptyset$.
\end{description}
\end{definition}

The following Darbo's fixed point theorem given in~\cite{bib1} is used in
the proof of the main result.

\begin{theorem}
\label{D} Let $M$ be a nonempty, bounded, convex and closed subset of the
space $E$ and let $T:M\rightarrow M$ be a continuous operator such that $\mu
(T(X))\leq k\mu (X)$ for all nonempty subset $X$ of $M$, where $k\in \lbrack
0,1)$ is a constant. Then $T$ has a fixed point in the subset $M$.
\end{theorem}

We consider the Banach space  $l^{\infty }$ of all real bounded sequences $%
x\colon {\mathbb{N}}_{k}\rightarrow {\mathbb{R}}$  equipped with the
standard supremum norm, i.e. 
\begin{equation*}
\Vert x\Vert =\sup_{n\in {\mathbb{N}}_{k}}|x_{n}|\text{ for }x\in \
l^{\infty }.
\end{equation*}%
Let $X$ be a nonempty, bounded subset of  $l^{\infty }$, $X_{n}=\left\{
x_{n}:x\in X\right\} $ (it means $X_{n}$ is a set of $n$-th terms of any
sequence belonging to $X$), and 
\begin{equation*}
diam\,\,X_{n}=\sup \left\{ \left\vert x_{n}-y_{n}\right\vert \colon x,y\in
X\right\} .
\end{equation*}%
We use a following measure of noncompactness in the space  $l^{\infty }$
(see~\cite{bibM1B}) 
\begin{equation*}
\mu (X)=\limsup_{n\rightarrow \infty }diam\,\,X_{n}.
\end{equation*}

\section{Main Result}

In this section, sufficient conditions for the existence of a bounded
solution of equation~\eqref{e0} are derived.

\begin{theorem}
\label{L2} Assume that $a,p,q:{\mathbb{N}}_{0}\rightarrow {\mathbb{R}}$, $r:{%
\mathbb{N}}_{0}\rightarrow {\mathbb{R}}\setminus \{0\}$, and $f:{\mathbb{R}}%
\rightarrow {\mathbb{R}}$. Let 
\begin{equation}  \label{alfa}
\alpha \geq 1 \,\, 
\mbox{is a ratio of positive integers with odd
denominator,}
\end{equation}
\begin{equation}  \label{gama}
\gamma \in \left( 0,1 \right] \,\, 
\mbox{ is a ratio of odd positive
integers,}
\end{equation}
and let $k$ be a fixed positive integer. Assume that 
\begin{equation}
f:{\mathbb{R}}\rightarrow {\mathbb{R}}\text{ is a locally Lipschitz function,%
}  \label{z1}
\end{equation}
and that the sequences $r:{\mathbb{N}}_{0}\rightarrow {\mathbb{R}}\setminus
\{0\}$, $a,q:{\mathbb{N}}_{0}\rightarrow {\mathbb{R}}$ satisfy 
\begin{equation}
\sum\limits_{n=0}^{\infty }\left\vert \frac{1}{r_{n}}\right\vert ^{\frac{1}{%
\gamma }}\sum\limits_{i=n}^{\infty }\left\vert a_{i}\right\vert <+\infty
\,\, \text{ and }\,\, \sum\limits_{n=0}^{\infty } \left\vert \frac{1}{r_{n}}%
\right\vert ^{\frac{1}{\gamma }}\sum\limits_{i=n}^{\infty }\left\vert
q_{i}\right\vert <+\infty .  \label{z2}
\end{equation}
Let the sequence $p:{\mathbb{N}}_{0}\rightarrow {\mathbb{R}}$ satisfies the
following condition 
\begin{equation}
-1<\liminf\limits_{n\rightarrow \infty }p_{n}\leq
\limsup\limits_{n\rightarrow \infty }p_{n}<1.  \label{z3}
\end{equation}
Assume additionally that 
\begin{equation}
\sum\limits_{i=0}^{\infty }\left\vert a_{i}\right\vert <+\infty, \,\,\,
\sum\limits_{i=0}^{\infty }\left\vert q_{i}\right\vert <+\infty.
\label{add_series}
\end{equation}
Then, there exists a bounded solution $x:{\mathbb{N}}_{k}\rightarrow {\ 
\mathbb{R}}$ of equation~\eqref{e0}.
\end{theorem}

Condition \eqref{z3} implies that there exist $n_{0}\in {\mathbb{N}}_{0}$
and a constant $P\in \lbrack 0,1)$ such that 
\begin{equation}
\left\vert p_{n}\right\vert \leq P<1,\text{ for }n\geq n_{0}.  \label{z5}
\end{equation}%
By condition~\eqref{add_series} there exists $n_{1}\in {\mathbb{N}}_{0}$
such that $\sum\limits_{i=n_{1}}^{\infty }\left\vert a_{i}\right\vert <1.$
Hence, by \eqref{gama}, 
\begin{equation*}
\sum\limits_{n=0}^{\infty }\left( \left\vert \frac{1}{r_{n}}\right\vert
\sum\limits_{i=n}^{\infty }\left\vert a_{i}\right\vert \right) ^{\frac{1}{%
\gamma }}\leq \sum\limits_{n=0}^{\infty }\left\vert \frac{1}{r_{n}}%
\right\vert ^{\frac{1}{\gamma }}\sum\limits_{i=n}^{\infty }\left\vert
a_{i}\right\vert .
\end{equation*}%
The above and condition \eqref{z2}, imply that 
\begin{equation}
\sum\limits_{n=0}^{\infty }\left( \left\vert \frac{1}{r_{n}}\right\vert
\sum\limits_{i=n}^{\infty }\left\vert a_{i}\right\vert \right) ^{\frac{1}{%
\gamma }}<+\infty .  \label{z31}
\end{equation}%
Analogously, we get 
\begin{equation}
\sum\limits_{n=0}^{\infty }\left( \left\vert \frac{1}{r_{n}}\right\vert
\sum\limits_{i=n}^{\infty }\left\vert q_{i}\right\vert \right) ^{\frac{1}{%
\gamma }}\leq \sum\limits_{n=0}^{\infty }\left\vert \frac{1}{r_{n}}%
\right\vert ^{\frac{1}{\gamma }}\sum\limits_{i=n}^{\infty }\left\vert
q_{i}\right\vert .  \label{z211}
\end{equation}%
Recalling that remainder of a series is the difference between the $n$--th
partial sum and the sum of a series, we denote by $\alpha _{n}$ and by $%
\beta _{n}$ the following remainders 
\begin{equation}
\alpha _{n}=\sum\limits_{j=n}^{\infty }\left( \left\vert \frac{1}{r_{j}}%
\right\vert \sum\limits_{i=j}^{\infty }\left\vert a_{i}\right\vert \right) ^{%
\frac{1}{\gamma }}\text{ and }\,\,\beta _{n}=\sum\limits_{j=n}^{\infty
}\left( \left\vert \frac{1}{r_{j}}\right\vert \sum\limits_{i=j}^{\infty
}\left\vert q_{i}\right\vert \right) ^{\frac{1}{\gamma }}  \label{ab}
\end{equation}%
We see, by \eqref{z31} and \eqref{z211} that 
\begin{equation}
\lim\limits_{n\rightarrow \infty }\alpha _{n}=0\,\,\text{ and }%
\lim\limits_{n\rightarrow \infty }\beta _{n}=0.  \label{z7}
\end{equation}%
Fix any number $d>0$. From \eqref{z1}, there exists a constant $M_{d}>0$
such that 
\begin{equation}
\left\vert f\left( x\right) \right\vert \leq M_{d}\mbox{  for  all }x\in %
\left[ -d,d\right] .  \label{Md}
\end{equation}%
Chose a constant $C$ such that 
\begin{equation}
0<C\leq \frac{d-Pd}{\left( 2^{\frac{1}{\gamma }-1}\left( M_{d}\right) ^{%
\frac{1}{\gamma }}+2^{\frac{1}{\gamma }-1}\left( d^{\alpha }\right) ^{\frac{1%
}{\gamma }}\right) }.  \label{C}
\end{equation}%
By condition~\eqref{z211} there exists a positive integer $n_{2}$ such that 
\begin{equation}
\alpha _{n}\leq C\text{ and }\beta _{n}\leq C\mbox{ for }n\in {\mathbb{N}}%
_{n_{2}}.  \label{z8}
\end{equation}%
\bigskip 

Define set $B$ as follows 
\begin{equation*}
B\colon =\left\{ (x_{n})_{n=0}^{\infty }:\left\vert x_{n}\right\vert \leq d,%
\text{ for }n\in {\mathbb{N}}_{n_{0}}\right\} ,
\end{equation*}%
Define a mapping $T\colon B\rightarrow l^{\infty }$ as follows 
\begin{equation}
(Tx)_{n}=%
\begin{cases}
-p_{n}x_{n-k}-\sum\limits_{j=n}^{\infty }\left( \frac{1}{r_{j}}%
\sum\limits_{i=j}^{\infty }\left( a_{i}f(x_{i+1})+q_{i}x_{i}^{\alpha
}\right) \right) ^{\frac{1}{\gamma }}, & \text{for any}\ n\geq n_{3}, \\ 
x_{n}, & \text{for any}\ 0\leq n<n_{3}%
\end{cases}
\label{z10}
\end{equation}%
where $n_{3}=\max \left\{ n_{1},n_{2}\right\} +k$. Observe that $B$ is a
nonempty, bounded, convex and closed subset of $l^{\infty }$. \newline
We will prove that the mapping $T$ has a fixed point in $B$. This proof will
follow in several subsequent steps.

\textbf{Step 1. Firstly, we show that }$T(B)\subset B$.

We will use classical inequality 
\begin{equation}  \label{ci}
\left( a+b\right) ^{s}\leq 2^{s-1}\left( a^{s}+b^{s}\right), \,\, a,b>0,
\,\, s\geq 1
\end{equation}
and the fact $t\rightarrow t^{1/\gamma }$ is nondecreasing. If $x\in B$,
then for $n< n_3$ $\left\vert (Tx)_{n}\right\vert=\left\vert
x_{n}\right\vert\leq d$ and by~\eqref{z10}, we get for any $n\geq n_3$ 
\begin{equation*}
\left\vert (Tx)_{n}\right\vert \,\,\leq \left\vert p_{n}\right\vert
\left\vert x_{n-k}\right\vert +\left\vert \sum\limits_{j=n}^{\infty }\left( 
\frac{1}{r_{j}}\sum\limits_{i=j}^{\infty }\left(
a_{i}f(x_{i+1})+q_{i}x_{i}^{\alpha }\right) \right) ^{\frac{1}{\gamma }
}\right\vert
\end{equation*}
\begin{equation*}
\leq \left\vert p_{n}\right\vert \left\vert x_{n-k}\right\vert
+\sum\limits_{j=n}^{\infty }\left( \left\vert \frac{1}{r_{j}}
\sum\limits_{i=j}^{\infty }\left( a_{i}f(x_{i+1})+q_{i}x_{i}^{\alpha
}\right) \right\vert \right) ^{\frac{1}{\gamma }}
\end{equation*}
\begin{equation*}
\leq \left\vert p_{n}\right\vert \left\vert x_{n-k}\right\vert
+\sum\limits_{j=n}^{\infty }\left( \left\vert \frac{1}{r_{j}}\right\vert
\sum\limits_{i=j}^{\infty }\left( \left\vert a_{i}\right\vert \left\vert
f(x_{i+1})\right\vert +\left\vert q_{i}\right\vert \left\vert
x_{i}\right\vert ^{\alpha }\right) \right) ^{\frac{1}{\gamma }}.
\end{equation*}
From \eqref{Md}, taking into account that $x_{n-k} \in B$, and because of $%
x_i \in B$ we have $\left\vert x_{i}\right\vert^\alpha \leq d^\alpha$. Thus 
\begin{equation*}
\left\vert (Tx)_{n}\right\vert \,\,\leq \left\vert p_{n}\right\vert d
+\sum\limits_{j=n}^{\infty }\left( \left\vert \frac{1}{r_{j}}\right\vert
\sum\limits_{i=j}^{\infty }\left( \left\vert a_{i}\right\vert
M_{d}+\left\vert q_{i}\right\vert d^{\alpha }\right) \right) ^{\frac{1}{%
\gamma }}
\end{equation*}
\begin{equation*}
\leq \left\vert p_{n}\right\vert d +\sum\limits_{j=n}^{\infty }\left(
\left\vert \frac{1}{r_{j}}\right\vert \left( \sum\limits_{i=j}^{\infty
}\left\vert a_{i}\right\vert M_{d}+\sum\limits_{i=j}^{\infty }\left\vert
q_{i}\right\vert d^{\alpha } \right) \right) ^{\frac{1}{\gamma }}.
\end{equation*}
By inequality \eqref{ci}, we have 
\begin{equation*}
\left\vert (Tx)_{n}\right\vert \,\,\leq \left\vert p_{n}\right\vert d +2^{%
\frac{1}{\gamma }-1}\sum\limits_{j=n}^{\infty }\left( \left( \left\vert 
\frac{1}{r_{j}}\right\vert \sum\limits_{i=j}^{\infty }\left\vert
a_{i}\right\vert M_{d}\right) ^{\frac{1}{\gamma }}+\left( \left\vert \frac{1%
}{r_{j}}\right\vert \sum\limits_{i=j}^{\infty }\left\vert q_{i}\right\vert
d^{\alpha }\right) ^{\frac{1}{\gamma }}\right)
\end{equation*}
\begin{equation*}
\leq\left\vert p_{n}\right\vert d+2^{\frac{1}{\gamma}-1}\left(M_{d}\right)^{%
\frac{1}{\gamma}}\sum\limits_{j=n}^{\infty}\left(\left\vert \frac{1}{r_{j}}%
\right\vert \sum\limits_{i=j}^{\infty}\left\vert a_{i}\right\vert \right)^{%
\frac{1}{\gamma}} + 2^{\frac{1}{\gamma}-1}\left( d^{\alpha}\right)^{\frac{1}{
\gamma}}\sum\limits_{j=n}^{\infty }\left( \left\vert \frac{1}{r_{j}}%
\right\vert \sum\limits_{i=j}^{\infty}\left\vert q_{i}\right\vert \right)^{%
\frac{1}{\gamma}}.
\end{equation*}
By using \eqref{z5}, \eqref{z8} and \eqref{C}, we estimate 
\begin{equation*}
\left\vert (Tx)_{n}\right\vert \,\,\leq Pd+2^{\frac{1}{\gamma }-1}\left(
M_{d}^{{}}\right) ^{\frac{1}{\gamma }}C+2^{\frac{1}{\gamma }-1}\left(
d^{\alpha }\right) ^{\frac{1}{\gamma }}C
\end{equation*}
\begin{equation*}
\leq Pd+\left( 2^{\frac{1}{\gamma }-1}\left( M_{d}^{{}}\right) ^{\frac{1}{%
\gamma }}+2^{\frac{1}{\gamma }-1}\left( d^{\alpha }\right) ^{\frac{1}{\gamma 
}}\right) \frac{d-Pd}{\left( 2^{\frac{1}{\gamma }-1}\left( M_{d}^{{}}\right)
^{\frac{1}{\gamma }}+2^{\frac{1}{\gamma }-1}\left( d^{\alpha }\right) ^{%
\frac{1}{\gamma }}\right) }=d.  \label{calc}
\end{equation*}
From above, we have estimation 
\begin{equation}  \label{d}
\left\vert (Tx)_{n}\right\vert \leq d, \mbox{ for }n\in {\mathbb{N}}_{n_{3}}.
\end{equation}

\textbf{Step 2. }$T$\textbf{\ is continuous}

By assumption \eqref{z1}, \eqref{add_series}, and by definition of $B$,
there exists a constant $c^*>0$ such that 
\begin{equation*}
\sum\limits_{i=j}^{\infty }\left\vert a_{i}f(x_{i+1})+q_{i}x_{i}^{\alpha
}\right\vert \leq c^*
\end{equation*}
for all $x\in B$. From \eqref{gama}, $t\rightarrow t^{1/\gamma }$ is locally
Lipschitz then it is Lipschitz on closed and bounded intervals. Hence, there
exists a constant $L_{\gamma }$ such that 
\begin{equation}  \label{Lgama}
\left\vert t^{1/\gamma }-s^{1/\gamma }\right\vert \leq L_{\gamma }\left\vert
t-s\right\vert \text{ for all }t,s\in \left[ -c^*,c^*\right].
\end{equation}
From \eqref{z1}, function $f$ is Lipschitz on $\left[ -d,d\right] $. So,
there is a constant $L_{d}>0$ such that 
\begin{equation}
\left\vert f\left( x\right) -f\left( y\right) \right\vert \leq
L_{d}\left\vert x-y\right\vert  \label{lipschitz}
\end{equation}
for all $x, y\in \left[ -d,d\right] $. From \eqref{alfa}, $x\rightarrow
x^{\alpha }$ is also Lipschitz on $\left[ -d,d\right] $. Then there is a
constant $L_{\alpha }$ such that 
\begin{equation}  \label{Lalfa}
\left\vert x^{\alpha }-y^{\alpha }\right\vert \leq L_{\alpha }\left\vert
x-y\right\vert \text{ for all }x,y\in \left[ -d,d\right].
\end{equation}

Let $(y^{(p)})$ be a sequence in $B$ such that $\left\Vert
y^{(p)}-x\right\Vert \rightarrow 0$ as $p\rightarrow \infty $. Since $B$ is
closed, $x\in B$. By \eqref{z10} and \eqref{calc}, we get 
\begin{equation*}
\begin{array}{l}
\forall\, n\geq n_3\quad \left\vert (Tx)_{n}-(Ty^{(p)})_{n}\right\vert
\,\,\leq \left\vert p_{n}\right\vert \left\vert
x_{n-k}-y_{n-k}^{(p)}\right\vert \\ 
\\ 
+\sum\limits_{j=n}^{\infty }\left\vert \frac{1}{r_{j}}\right\vert ^{\frac{1}{
\gamma }}\left\vert \left(\sum\limits_{i=j}^{\infty }\left(
a_{i}f(x_{i+1})+q_{i}\left( x_{i}\right) ^{\alpha }\right) \right) ^{\frac{1%
}{\gamma }}-\left( \sum\limits_{i=j}^{\infty }\left(
a_{i}f(y_{i+1}^{(p)})+q_{i}\left( y_{i}^{\left( p\right) }\right) ^{\alpha
}\right) \right) ^{\frac{1}{\gamma }}\right\vert.%
\end{array}%
\end{equation*}
From \eqref{Lgama}, we have 
\begin{equation*}
\begin{array}{l}
\forall\, n\geq n_3\quad \left\vert (Tx)_{n}-(Ty^{(p)})_{n}\right\vert
\,\,\leq \bigskip \left\vert p_{n}\right\vert \left\vert
x_{n-k}-y_{n-k}^{(p)}\right\vert \\ 
\\ 
+\sum\limits_{j=n}^{\infty }\left\vert \frac{1}{r_{j}}\right\vert ^{\frac{1}{%
\gamma }}L_{\gamma }\left\vert \sum\limits_{i=j}^{\infty
}a_{i}f(x_{i+1})+\sum\limits_{i=j}^{\infty }q_{i}\left( x_{i}\right)
^{\alpha }-\sum\limits_{i=j}^{\infty }a_{i}f(y_{i+1}^{\left( p\right)
})-\sum\limits_{i=j}^{\infty }q_{i}\left( y_{i}^{\left( p\right) }\right)
^{\alpha }\right\vert \\ 
\\ 
\leq \bigskip \left\vert p_{n}\right\vert \left\vert
x_{n-k}-y_{n-k}^{(p)}\right\vert +L_{\gamma }\sum\limits_{j=n}^{\infty
}\left\vert \frac{1}{r_{j}}\right\vert ^{\frac{1}{\gamma }%
}\sum\limits_{i=j}^{\infty }\left\vert a_{i}\right\vert \left\vert
f(x_{i+1})-f(y_{i+1}^{\left( p\right) })\right\vert \\ 
\\ 
+L_{\gamma }\sum\limits_{j=n}^{\infty }\left\vert \frac{1}{r_{j}}\right\vert
^{\frac{1}{\gamma }}\sum\limits_{i=j}^{\infty }\left\vert q_{i}\right\vert
\left\vert \left( x_{i}\right) ^{\alpha }-\left( y_{i}^{\left( p\right)
}\right) ^{\alpha }\right\vert.%
\end{array}%
\end{equation*}
Hence, by \eqref{lipschitz} and \eqref{Lalfa}, we obtain 
\begin{equation*}
\begin{array}{l}
\forall\, n\geq n_3\quad \left\vert (Tx)_{n}-(Ty^{(p)})_{n}\right\vert
\,\,\leq \left\vert p_{n}\right\vert \left\vert
x_{n-k}-y_{n-k}^{(p)}\right\vert \bigskip \\ 
+L_{\gamma }L_{d}\sum\limits_{j=n}^{\infty }\left\vert \frac{1}{r_{j}}%
\right\vert ^{\frac{1}{\gamma }}\sum\limits_{i=j}^{\infty }\left\vert
a_{i}\right\vert \left\vert x_{i+1}-y_{i+1}^{\left( p\right) }\right\vert
+L_{\gamma }L_{\alpha }\sum\limits_{j=n}^{\infty }\left\vert \frac{1}{r_{j}}%
\right\vert ^{\frac{1}{\gamma }}\sum\limits_{i=j}^{\infty }\left\vert
q_{i}\right\vert \left\vert x_{i}-y_{i}^{\left( p\right) }\right\vert
\bigskip \\ 
\leq \sup\limits_{i \in \mathbb{N}_0}\left\vert y^{(p)}_i-x_i\right\vert
\left( \left\vert p_{n}\right\vert +L_{\gamma
}L_{d}\sum\limits_{j=n}^{\infty }\left\vert \frac{1}{r_{j}}\right\vert ^{%
\frac{1}{\gamma }}\sum\limits_{i=j}^{\infty }\left\vert a_{i}\right\vert
+L_{\gamma }L_{\alpha }\sum\limits_{j=n}^{\infty }\left\vert \frac{1}{r_{j}}%
\right\vert ^{\frac{1}{\gamma }}\sum\limits_{i=j}^{\infty }\left\vert
q_{i}\right\vert \right).%
\end{array}%
\end{equation*}
Moreover, 
\begin{equation*}
\forall\, 0\leq n< n_3\quad \left\vert (Tx)_{n}-(Ty^{(p)})_{n}\right\vert
\,\,\leq \left\Vert y^{(p)}-x\right\Vert
\end{equation*}
Thus, by \eqref{z2} and \eqref{z3}, we have 
\begin{equation*}
\lim\limits_{p\rightarrow \infty }\left\Vert Ty^{(p)}-Tx\right\Vert =0\text{
as }\lim\limits_{p\rightarrow \infty }\left\Vert y^{(p)}-x\right\Vert =0.
\end{equation*}
This means that $T$ is continuous.

\textbf{Step 3. Comparison of the measure of noncompactness}

Now, we need to compare a measure of noncompactness of any subset $X$ of $B$
and $T(X)$. Let us fix any nonempty set $X\subset B$. Take any sequences $%
x,y\in X$. Following the same calculations which led to the continuity of
the operator $T$ we see that  
\begin{equation*}
\forall \,n\geq n_{3}\quad \left\vert (Tx)_{n}-(Ty)_{n}\right\vert \leq
\left\vert p_{n}\right\vert \left\vert x_{n-k}-y_{n-k}\right\vert +L_{\gamma
}L_{\alpha }\beta _{n}\left\vert x_{n}-y_{n}\right\vert +L_{\gamma
}L_{d}\alpha _{n}\left\vert x_{n+1}-y_{n+1}\right\vert .
\end{equation*}%
Taking sufficiently large $n$, by \eqref{ab} and \eqref{z7}, we get 
\begin{equation*}
L_{\gamma }L_{d}\alpha _{n}\leq c_{1}<\frac{1-P}{4}\,\,\,\mbox{ and }%
\,\,\,L_{\gamma }L_{\alpha }\beta _{n}\leq c_{2}<\frac{1-P}{4}
\end{equation*}%
Here $c_{1},\,\,c_{2}$ are some real constants. From \eqref{z5}, we have 
\begin{equation*}
P+c_{1}+c_{2}<\frac{1+P}{2}.
\end{equation*}%
We see that exists $n_{5}$ such that 
\begin{equation*}
\forall \,n\geq n_{5}\quad diam\,\,(T(X))_{n}\leq
Pdiam\,\,X_{n-k}+c_{1}diam\,\,X_{n}+c_{2}diam\,\,X_{n+1}.
\end{equation*}%
This yields by the properties of the upper limit that 
\begin{equation*}
\limsup_{n\rightarrow \infty }diam\,\,(T(X))_{n}\leq
\,P\,\limsup_{n\rightarrow \infty
}diam\,\,X_{n-k}+c_{1}\,\limsup_{n\rightarrow \infty
}diam\,\,X_{n}+c_{2}\limsup_{n\rightarrow \infty }diam\,\,X_{n+1}.
\end{equation*}%
From above, for any $X\subset B$, we have $\mu (T(X))\leq \left(
c_{1}\,+c_{2}+P\right) \mu (X)$.

\textbf{Step 4. Relation between fixed points and solutions}

By Theorem~\ref{D} we conclude that $T$ has a fixed point in the set $B$. It
means that there exists $x\in B$ such that 
\begin{equation*}
x_{n}=(Tx)_{n}.
\end{equation*}%
Thus 
\begin{equation}
x_{n}=-p_{n}x_{n-k}-\sum\limits_{j=n}^{\infty }\left( \frac{1}{r_{j}}%
\sum\limits_{i=j}^{\infty }\left( a_{i}f(x_{i+1})+q_{i}x_{i}^{\alpha
}\right) \right) ^{\frac{1}{\gamma }},\,\text{ for }\,\,n\in {\mathbb{N}}%
_{n_{3}}  \label{z12}
\end{equation}

To show that there exists a correspondence between fixed points of $T$ and
solutions to ~\eqref{e0} we apply operator $\Delta $ to both sides of the
following equation 
\begin{equation*}
x_{n}+p_{n}x_{n-k}=-\sum\limits_{j=n}^{\infty }\left( \frac{1}{r_{j}}%
\sum\limits_{i=j}^{\infty }\left( a_{i}f(x_{i+1})+q_{i}x_{i}^{\alpha
}\right) \right) ^{\frac{1}{\gamma }},
\end{equation*}%
which is obtained from~\eqref{z12}. We find that 
\begin{equation*}
\Delta (x_{n}+p_{n}x_{n-k})=\left( \frac{1}{r_{n}}\sum\limits_{i=n}^{\infty
}\left( a_{i}f(x_{i+1})+q_{i}x_{i}^{\alpha }\right) \right) ^{\frac{1}{%
\gamma }},\,\,\,n\in {\mathbb{N}}_{n_{3}}.
\end{equation*}%
and next 
\begin{equation*}
\left( \Delta (x_{n}+p_{n}x_{n-k})\right) ^{\gamma }=\frac{1}{r_{n}}%
\sum\limits_{i=n}^{\infty }\left( a_{i}f(x_{i+1})+q_{i}x_{i}^{\alpha
}\right) ,\,\,\,n\in {\mathbb{N}}_{n_{3}}.
\end{equation*}%
Taking operator $\Delta $ again to both sides of the above equation we
obtain 
\begin{equation*}
\Delta \left( r_{n}\left( \Delta (x_{n}+p_{n}x_{n-k})\right) ^{\gamma
}\right) =-a_{n}f(x_{n+1})-q_{n}x_{n}^{\alpha },\,\,\,n\in {\mathbb{N}}%
_{n_{3}}.
\end{equation*}%
So, we get equation~\eqref{e0} for $n\in {\mathbb{N}}_{n_{3}}$. Sequence $x$%
, which is a fixed point of mapping $T$, is a bounded sequence which
fulfills equation~\eqref{e0} for large $n$. If $n_{3}\geq k$ the proof is
ended. We find previous $n_{3}-k+1$ terms of sequence $x$ by formula 
\begin{equation*}
x_{n-k+l}=\frac{1}{p_{n+l}}\left( -x_{n+l}+\sum\limits_{j={n+l}}^{\infty
}\left( \frac{1}{r_{j}}\sum\limits_{i=j}^{\infty }\left(
a_{i}f(x_{i+1})+q_{i}x_{i}^{\alpha }\right) \right) ^{\frac{1}{\gamma }%
}\right) ,
\end{equation*}%
where $l\in \left\{ 0,1,2,\dots ,k-1\right\} $, which results leads directly
from~\eqref{e0}. It means that equation~\eqref{e0} has at least one bounded
solution $x:{\mathbb{N}}_{k}\rightarrow {\mathbb{R}}$.

This completes the proof.

\textbf{Remark. }\textit{We note the previous terms of the solution sequence
are not obtained through a fixed point method, but through backward
iteration. It is common that one has a }$1-1$\textit{\ correspondence
between fixed points to a suitably chosen operator and solutions to the
problem under consideration. Here we get as a fixed point solution some
sequence which starting from some index is a solution to the given problem
and in which the first terms must be iterated. This procedure must be
applied since we see that in equation~\eqref{e0} we have to know also
earlier terms in order to start iteration; this is the so called iteration
with memory. We recall that in recent works concerning application of the
measure of noncompactness to discrete equations, only problems without
memory have been considered. That is why we had to alter to established
procedure to overcome the difficulty arising in this problem. We believe our
method would be applicable for several other problems}

\section{A special type stability}

The type of stability investigated in this paper is contained in the
following theorem.

\begin{theorem}
\label{T1} Assume that 
\begin{equation}
q_{n}\equiv 0,  \label{e}
\end{equation}%
and conditions \eqref{gama} and \eqref{z2}--\eqref{add_series} are held.
Assume further that there exists a positive constant $D$ such that 
\begin{equation*}
\left\vert f(u)-f(v)\right\vert \leq D\left\vert u-v\right\vert 
\end{equation*}%
for any $u,v\in {\mathbb{R}}$. Then equation~\eqref{e0} has at least one
solution $x:{\mathbb{N}}_{k}\rightarrow {\mathbb{R}}$ with the following
stability property: given any other solution $y:{\mathbb{N}}_{k}\rightarrow {%
\mathbb{R}}$ and $\varepsilon >0$ there exists $T>$ $n_{3}$ such that for
every $t\geq T$ the following inequality holds 
\begin{equation*}
\left\vert x(t)-y(t)\right\vert \leq \varepsilon .
\end{equation*}
\end{theorem}

From Theorem~\ref{L2}, equation~\eqref{e0} has at least one bounded solution 
$x:{\mathbb{N}}_{0}\rightarrow {\mathbb{R}}$ which can be rewritten in the
form 
\begin{equation*}
x_{n}=(Tx)_{n},
\end{equation*}%
where mapping $T$ is defined by~\eqref{z10} for $n\geq n_{3}$. By Definition %
\ref{def2}, sequence $x$ is an asymptotically stable solution of equation $%
x_{n}=(Tx)_{n}$ From the above and condition \eqref{e}, analogously as the
steps in the proof of Theorem \ref{L2}, we see that 
\begin{equation*}
\begin{array}{l}
\left\vert x_{n}-y_{n}\right\vert =\left\vert (Tx)_{n}-(Ty)_{n}\right\vert
\,\leq \bigskip \\ 
\left\vert p_{n}\right\vert \left\vert x_{n-k}-y_{n-k}\right\vert +L_{\gamma
}D\sum\limits_{j=n}^{\infty }\left\vert \frac{1}{r_{j}}\right\vert ^{\frac{1%
}{\gamma }}\sum\limits_{i=j}^{\infty }\left\vert a_{i}\right\vert \left\vert
x_{i+1}-y_{i+1}\right\vert .%
\end{array}%
\end{equation*}%
Note that for $n$ large enough, say $n\geq n_{4}\geq n_{3}$, we have 
\begin{equation*}
\vartheta :=\left\vert p_{n}\right\vert +L_{\gamma
}D\sum\limits_{j=n}^{\infty }\left\vert \frac{1}{r_{j}}\right\vert ^{\frac{1%
}{\gamma }}\sum\limits_{i=j}^{\infty }\left\vert a_{i}\right\vert <1
\end{equation*}%
Let us denote 
\begin{equation*}
\limsup_{n\rightarrow \infty }\left\vert x_{n}-y_{n}\right\vert =l,
\end{equation*}%
and observe that 
\begin{equation*}
\limsup_{n\rightarrow \infty }\left\vert x_{n}-y_{n}\right\vert
=\limsup_{n\rightarrow \infty }\left\vert x_{n-k}-y_{n-k}\right\vert
=\limsup_{n\rightarrow \infty }\left\vert x_{n+1}-y_{n+1}\right\vert .
\end{equation*}%
Thus, from the above, we have 
\begin{equation*}
l\leq \vartheta \cdot l.
\end{equation*}%
This means that $\limsup\limits_{n\rightarrow \infty }\left\vert
x_{n}-y_{n}\right\vert =0$. This completes the proof since for $\varepsilon
>0$ there exists $n_{4}\in {\mathbb{N}}_{0}$ such that for every $n\geq
n_{4}\geq n_{3}$ the following inequality holds 
\begin{equation*}
\left\vert x_{n}-y_{n}\right\vert \leq \varepsilon .
\end{equation*}

\section{Comments and an example}

In \cite{SchmeildelZbaszyniakCAMW}, the authors consider a special type of
problem \eqref{e0}, namely they investigate the existence of a solution and
Lyapunov type stability to the following equation 
\begin{equation}
\Delta \left( r_{n}\Delta x_{n}\right) =a_{n}f(x_{n+1}).  \label{ESZZ_RPW}
\end{equation}
Their main assumption is the linear growth assumption on nonlinear term $f$.
More precisely, they assume that there exists a positive constant $M$ such
that $\left\vert f(x_{n})\right\vert \leq M$ $\left\vert x_{n}\right\vert $
for all $x\in N_{0}$. Using ideas developed in this paper we get the
following result.

\begin{theorem}
Assume that $f:{\mathbb{R}}\rightarrow {\mathbb{R}}$ satisfies the condition %
\eqref{z1} and the sequences $r:{\mathbb{N}}_{0}\rightarrow {\mathbb{R}}%
\setminus \{0\}$, $a:{\mathbb{N}}_{0}\rightarrow {\mathbb{R}}$ are such that 
\begin{equation*}
\sum\limits_{n=0}^{\infty }\left\vert \frac{1}{r_{j}}\right\vert
\sum\limits_{i=n}^{\infty }\left\vert a_{i}\right\vert <+\infty.
\end{equation*}
Then, there exists a bounded solution $x:{\mathbb{N}}_{0}\rightarrow {%
\mathbb{R}}$ of equation~\eqref{ESZZ_RPW}.
\end{theorem}

Finally, we give an example of equation which can be considered by our
method.

\begin{example}
Take $k=3$, an arbitrary $C^{1}$ function $f:{\mathbb{R\rightarrow R}}$ and
consider the following problem 
\begin{equation}
\Delta \left( \left( -1\right) ^{n}\Delta \left( x_{n}+\frac{1}{2}%
x_{n-3}\right) ^{1/3}\right) +\frac{1}{2^{n}} \left( \left( x_{n}\right)
^{5}+f\left( x_{n+1}\right)\right) =0.  \label{prob}
\end{equation}
Taking $\gamma =\frac{1}{3}$, $\alpha =5$, $r_{n}=\left( -1\right) ^{n}$, $%
p_{n}=\frac{1}{2}$, $a_{n}=q_{n}=\frac{1}{2^{n}}$ with $f\left( x\right)
=x^5 $ we see that $x_{n}=\left( -1\right) ^{n}$ is a bounded solution to %
\eqref{prob}. By Theorem \ref{T1}, this solution is asymptotically stable.
\end{example}

\end{document}